\documentclass[aps,pre,twocolumn,groupedaddress,amsmath,amssymb]{revtex4}

\usepackage{graphicx}

\begin{document}


\title{The Type II Phase Resetting Curve is Optimal for Stochastic Synchrony}


\author{Aushra Abouzeid and Bard Ermentrout}
\affiliation{University of Pittsburgh}


\date{\today}

\begin{abstract}
The phase-resetting curve (PRC) describes the response of a neural
oscillator to small perturbations in membrane potential. Its
usefulness for predicting the dynamics of weakly coupled deterministic
networks has been well characterized. However, the inputs to real
neurons may often be more accurately described as barrages of synaptic
noise. Effective connectivity between cells may thus arise in the form
of correlations between the noisy input streams. We use constrained
optimization and perturbation methods to prove that PRC shape
determines susceptibility to synchrony among otherwise uncoupled
noise-driven neural oscillators. PRCs can be placed into two general
categories: Type I PRCs are non-negative while Type II PRCs have a
large negative region. Here we show that oscillators with
Type~II PRCs receiving common noisy input sychronize more readily than
those with Type~I PRCs. 
\end{abstract}


\maketitle


\section*{Introduction}

Synchronous oscillations are found in many brain areas and are
responsible for macroscopic electrical responses of the brain
including field potentials and EEG signals. Within a single brain
area, synchronization of neuronal activity serves to amplify signals
to upstream regions \cite{tiesinga:2004}, while synchronization
across different areas may allow activity to be selectively routed. 

Considerable theoretical interest has recently emerged in the
generation of synchrony by correlated ``noisy'' inputs to uncoupled
oscillators \cite{teramae:2004, goldobin:2005, nakao:2005, stroeve:2001},
a phenomenon we will refer to as stochastic synchrony. In the brain,
stochastic synchrony may account for observations such as long-range
synchronization \cite{engel:1991, engel:1991b}, that are
difficult to explain by the presence of synaptic connectivity
alone. Moreover, noisy inputs have been shown to synchronize real
neurons \textit{in vitro} \cite{galan:2006}.

The key component in the study of noisy oscillators is the
phase-resetting curve (PRC). This curve characterizes how inputs to an
oscillator shift its timing, or phase. In the context of neurons, spike
times are believed to play an important role in coding and in the propagation of information across brain regions. Thus, the PRC
provides a quantitative characterization of how inputs to neural
oscillators alter the timing of spikes.

The theory of deterministic oscillators has shown that
the type of bifurcation from steady-state to periodic behavior
determines the shape of the PRC. Weak coupling theory shows that the form of the
interaction between oscillators together with their intrinsic response (the PRC)
provide sufficient information about the ability of the coupling to
synchronize (or desynchronize) the oscillations. For very fast
excitatory synaptic interactions, Type II oscillators characterized
by the Hopf bifurcation synchronize more readily than Type I 
oscillators characterized by the
saddle-node-on-an-invariant-circle (SNIC) bifurcation
\cite{hansel:1995, ermentrout:2001, gutkin:2005, netoff:2005}. This difference 
in ability to synchronize with excitatory
coupling is a consequence of the shape of the PRC occurring near the
two different bifurcations. A PRC which contains both negative and
positive lobes can allow inputs to both slow down the oscillator which is ahead and
speed up the oscillator which is behind. In contrast, a non-negative
PRC can only speed up the timing of both oscillators, so that
synchronization becomes more difficult. A number of authors
\cite{ermentrout:2001, izhikevich, brown} have shown that the PRC near
a SNIC is non-negative and approximately proportional to $1-\cos t$,
while the PRC near a Hopf is proportional to $\sin(t+\alpha).$ Thus, Type II 
PRCs have a large negative lobe, whereas Type I PRCs are strictly positive. 

Two recent papers have shown that Type II PRCs
are better than Type I PRCs 
at synchronizing uncoupled oscillators with correlated input \cite{galanfp,sashi}. That is,
for a given input correlation of the noisy stimulus, the output
correlation of the oscillators is higher with Type II than with Type I
PRCs. In these two papers, specific functions for PRCs were checked
(namely, $\sin(t)$ and $1-\cos(t)$), and the correlations and degree of
synchrony were analytically and numerically computed. However, it is
not known whether there are other PRC shapes that might produce even stronger
stochastic synchronization. 

The easiest way to quantify stochastic synchrony 
is to examine the Lyapunov exponent, the rate at which two oscillators
receiving identical inputs converge to synchrony. In this paper we will explore how
this quantity depends on the shape of the PRC. In particular, we find
that Type II PRCs lead to faster convergence than do Type I, and
we use variational principles to determinine the optimal shape of the
PRC to maximize this convergence. 

First in Section \ref{sec:ito} we introduce the phase reduction of a
stochastically driven neural oscillator using the It\^o change of
variables, and in Section \ref{sec:lyapunov} we derive the Lyapunov
exponent for two such oscillators receiving common noise. Next we use
the Fokker-Planck equation in Section \ref{sec:phase} to obtain the
probability distribution of the phase of a noise-driven neural
oscillator. The Euler-Lagrange method for constrained optimization
allows us in Section \ref{sec:euler-lagrange} to find the PRC that
minimizes the Lyapunov exponent. This leads to a 4th order system of
nonlinear differential equations, which we approximate to an arbitrary
order of accuracy using regular perturbations in Section
\ref{sec:perturbation}. The resulting approximation shows that a Type
II PRC achieves the minimal Lyapunov exponent, hence producing more
robust convergence to synchrony than a Type I PRC. Several interesting
cases that arise as a function of the constraint parameters are
discussed in Section \ref{sec:constraints}. Finally in Section
\ref{sec:numerics} we show that numerical solution of the 4th order
system agrees with the perturbation-derived approximation. 

\section{\label{sec:ito} It\^o Phase Reduction}

Consider a neural oscillator with additive white noise decribed by the stochastic differential equation
\begin{equation}
dX = F(X) dt + \sigma M dW,
\label{eqn:SDE}
\end{equation}
where $F(X)$ represents the deterministic equations of motion, 
$\sigma$ is the amplitude of the noise, $M$ is a constant matrix, and
$dW$ is a vector of Gaussian 
white noise. Note that for a general limit-cycle oscillator,
there need be no constraints on the entries of $M$. For neural models however,
the noise typically occurs in current felt by the neuron, and this current
appears only in the voltage-component of the deterministic
model. Without loss of generality, we take the voltage to be the first
component. Thus, we will assume here that $M$ has all zero entries except
for the $(1,1)$ element, which is identically 1. 

The phase reduction method \cite{teramae:2004} applied to
Eq.(\ref{eqn:SDE}) gives a stochastic differential equation for the
evolution of the oscillator's phase: 
\begin{equation}
d\theta = dt + \sigma \Delta(\theta) dW,
\label{eqn:Stratonovich}
\end{equation}
where we have assumed without loss of generality that the intrinsic
frequency of the oscillator is $\omega = 1$, and $dW$ is now a scalar
white noise process. Here $\Delta$ is the
infinitesimal phase response curve defined by 
\[
\Delta(\theta) := \nabla_X \theta \Big\vert_{X=X_0(\theta)}
\]
where $X_0(\theta)$ is the unperturbed limit-cycle solution of the
deterministic equation $\dot{X} = F(X)$. See Kuramoto \cite{kuramoto:1984}, pages 26-27.

It is now important to note that the usual phase reduction method uses
the conventional change of variables, so Eq.(\ref{eqn:Stratonovich})
must be regarded as a Stratonovich differential equation
\cite{teramae:2004, horsthemke:1984}. To eliminate the
correlation between $\theta$ and the white noise $\xi=dW$, we must
apply It\^{o}'s Lemma to obtain an equivalent but analytically more
convenient formulation 
\begin{equation}
d\theta = \left[ 1 + \frac{\sigma ^2}{2} \Delta'(\theta) \Delta(\theta)\right]dt + \sigma \Delta(\theta) dW,
\label{eqn:Ito}
\end{equation}
where $'$ denotes $\frac{\partial}{\partial \theta}$.
In a recent paper, Yoshimura and Arai \cite{yoshi} show that Eq.(\ref{eqn:Ito}) is incomplete and that another term must be added in
the case where the noise is strictly white. However, more recently (in
preparation) we show that the correct reduction is more subtle, and
under some reasonable circumstances the additional term can be made
arbitrarily small. Thus we will stay with the conventional
phase-reduced model as first proposed by Teramae and Tanaka \cite{teramae:2004}.  

\section{\label{sec:lyapunov} Lyapunov Exponent}

As a standard measure of susceptibility to synchrony, we will now
derive the Lyapunov exponent for two identical uncoupled neural
oscillators receiving common additive white noise. The resulting
analysis, however, applies equally well to an arbitrary number of
identical noninteracting oscillators. 

Let us define the phase difference $\phi := \theta_2-\theta_1$, where
$\theta_1$ and $\theta_2$ each obey Eq.(\ref{eqn:Ito}). Linearizing
around the synchronous state $\phi=0$, we obtain as in \cite{teramae:2004}:
\[
d\phi = \frac{\sigma ^2}{2} \left[(\Delta'\Delta)'(\theta)\phi\right]dt + \sigma [ \Delta'(\theta)\phi] dW,
\]
where $\theta$ obeys Eq.(\ref{eqn:Ito}) as well. Since the Lyapunov exponent
is defined as $\lambda := \lim_{t\rightarrow \infty}
\frac{\log(\phi(t))}{t}$, let us make the change of variables
$y:=\log(\phi)$. Once again we invoke It\^{o}'s Lemma, and after
simplification we find that $y$ satisfies the stochastic differential
equation 
\[
dy = \frac{\sigma^2}{2}[\Delta''\Delta] dt + \sigma \Delta' dW.
\]
Next we integrate, divide by $t$ and take the limit as $t\to\infty$ 
to obtain an expression for $\lambda$.
\begin{eqnarray*}
\lambda &=& \lim_{t\rightarrow \infty}\frac{y(t)}{t} \\
    &=& \lim_{t\rightarrow \infty} \frac{\sigma^2}{2t}\int_0^t \Delta''(\theta(s))\Delta(\theta(s))ds \\
&& + \frac{\sigma}{t} \int_0^t \Delta'(\theta(s)) dW(s)
\end{eqnarray*}
Assuming the system is ergodic, we can replace the long time average
on the right hand side with the spatial or ensemble average. Due to
the It\^{o} change of variables, the last term drops out leaving 
\begin{equation}
\lambda = \frac{\sigma^2}{2}\int_0^1 \Delta''(\theta) \Delta(\theta) P(\theta) d\theta,
\label{eqn:lambda}
\end{equation}
where $P(\theta)$ is the steady-state distribution of the phase. 

Note that Teramae and Tanaka derive an expression for $\lambda$ in \cite{teramae:2004} by making the approximation $P(\theta)=1$. Substituting this value into Eq.(\ref{eqn:lambda}) and performing integration by parts, they obtain

\[
\lambda \approx - \frac{\sigma^2}{2}\int_0^1(\Delta'(\theta))^2 d\theta.
\]
In this paper, however, we wish to retain the generality of $P(\theta)$ as discussed below.

\section{\label{sec:phase} Steady-State Phase Distribution} 

In order to evaluate the Lyapunov exponent, we need to
obtain the stationary density of the phase when perturbed by
noise. Teramae and Tanaka \cite{teramae:2004} have treated the density as uniform,
which is correct for weak noise. However our subsequent perturbation analysis will require higher-order terms, so we will need to derive a more accurate value for the steady-state phase distribution.

By applying the Fokker-Planck equation to (\ref{eqn:Ito}), we obtain
after simplification a partial differential equation for the
probability distribution $P(\theta,t)$: 
\[
\frac{\partial P}{\partial t} = -\frac{\partial P}{\partial \theta} + \frac{\sigma^2}{2}\frac{\partial }{\partial \theta} \left[\Delta \frac{\partial (\Delta P)}{\partial \theta}\right].
\]
Now we may set $\frac{\partial P}{\partial t} = 0$ to find the steady state, then integrate once with respect to $\theta$ to obtain:
\begin{equation}
-J = - P + \frac{\sigma^2}{2} \left[\Delta \frac{\partial (\Delta P)}{\partial \theta}\right],
\label{eqn:steadystate}
\end{equation}
where $-J$ is a constant of integration. We require that $P(0)=P(1)$
and that the solution be normalized, namely $\int_0^1P(\theta)\ d\theta=1.$ Note that the equations are singular, since $\Delta(\theta)$ generally vanishes at several
places, in particular at $\theta=0, 1$. In the appendix below, we prove the 
existence of the stationary density by directly solving the linear equations 
and taking appropriate limits.

In the remainder of this section, we use regular perturbation theory to approximate the
stationary density for small noise, $0<\sigma\ll 1.$ 
 To approximate both $J$ and $P$ we substitute
\begin{eqnarray*}
J &=& 1 + \sigma^2 J_1 + \sigma^4J_2 + \cdots \\
P(\theta) &=& 1 + \sigma^2 P_1(\theta) + \sigma^4 P_2(\theta) + \cdots
\end{eqnarray*}
into equation (\ref{eqn:steadystate}).
Equating like powers of $\sigma$ gives
\[
-J_1 = -P_1(\theta) + \frac{1}{2}\Delta(\theta) \Delta'(\theta).
\]
Integrating both sides over $[0,1]$ leaves the constant on the
left hand side unchanged. For the right hand side, 
note that $\int_0^1 P(\theta) d\theta = 1$, and 
hence $\int_0^1 P_1(\theta) d\theta = 0$. Furthermore, 
$\Delta\Delta' = \frac{1}{2}\frac{d}{d\theta}(\Delta^2)$ so that
\begin{eqnarray*}
J_1 &=& -\frac{1}{4}(\Delta(1)^2-\Delta(0)^2) \\
  &=& 0,
\end{eqnarray*}
since $\Delta$ is periodic. Thus we have $P_1(\theta) = \frac{1}{2}\Delta(\theta) \Delta'(\theta)$.

Similarly, 
\[
-J_2 = - P_2(\theta) + \frac{1}{2}\Delta(\theta)^2\Delta'(\theta)^2 + \frac{1}{4}\Delta(\theta)^3\Delta''(\theta).
\]
Since $\int_0^1 P_2(\theta) d\theta = 0$ as well, we can integrate both sides as above and use integration by parts to obtain
\begin{eqnarray*}
J_2 &=& \frac{1}{4}\int_0^1(\Delta(\theta)\Delta'(\theta))^2 d\theta \\
P_2(\theta) &=& \frac{1}{2}\Delta(\theta)^2\Delta'(\theta)^2 + \frac{1}{4}\Delta(\theta)^3\Delta''(\theta) \\
&& + \frac{1}{4}\int_0^1(\Delta(\theta)\Delta'(\theta))^2 d\theta.
\end{eqnarray*}

In summary,
\begin{eqnarray}
\nonumber
J &=& 1 + \frac{\sigma^4}{4}\int_0^1(\Delta(\theta)\Delta'(\theta))^2 d\theta \\
\nonumber
P(\theta) &=& 1 + \frac{\sigma^2}{2}\Delta(\theta) \Delta'(\theta) + \frac{\sigma^4}{4} \Bigg[ 2\Delta(\theta)^2\Delta'(\theta)^2 \\
&& + \Delta(\theta)^3\Delta''(\theta) + \int_0^1(\Delta(\theta)\Delta'(\theta))^2 d\theta \Bigg].
\label{eqn:ptosigma^4}
\end{eqnarray}
For the perturbation expansions in the next section, it will suffice to write $J = 1$. We will use Eq.(\ref{eqn:ptosigma^4}) in Section \ref{sec:constraints} and for the numerical verifications in Section \ref{sec:numerics}.

\section{\label{sec:euler-lagrange} Constrained Optimization}

The Euler-Lagrange variational technique provides a method for
determining the phase resetting curve $\Delta$ that minimizes the
Lyapunov exponent, subject to appropriate constraints. To ensure
smooth solutions and to eliminate uninformative harmonics of the
optimal solution, we begin by imposing the general constraint 
\begin{equation}
\int_0^1 a(\Delta(\theta))^2 + b(\Delta'(\theta))^2 +c(\Delta''(\theta))^2 d\theta = 1,
\label{eqn:constraint}
\end{equation} 
where $a$, $b$ and $c$ are free parameters. A standard normalization
has $a=1,b=0,c=0$, but non-zero values of $b$, $c$ endow solutions with
additional smoothness. Below we will explore the
cases that arise from specific choices of these. 

We proceed by placing Eqs.(\ref{eqn:lambda}), (\ref{eqn:steadystate}) and (\ref{eqn:constraint}) together with the approximation $J=1$ into the Euler-Lagrange formula to obtain the functional
\begin{eqnarray}
&\int_0^1 &\Delta'' \Delta P + \nu_1 \left[ a\Delta^2 + b(\Delta')^2 +c(\Delta'')^2 -1\right] \nonumber \\
&& + \nu_2(\theta) \left[1 -P + \frac{\sigma^2}{2} \Delta (\Delta P)' \right] d\theta = 0,
\label{eqn:EL-functional}
\end{eqnarray}
where $\nu_1$ is a free parameter, and $\nu_2(\theta)$ represents a continuum of free parameters. 

Define the operator
\begin{eqnarray*}
\mathcal{L}(\Delta) &:=& \Delta'' \Delta P + \nu_1 \left[ a\Delta^2 + b(\Delta')^2 +c(\Delta'')^2-1 \right] \\
&& + \nu_2(\theta) \left[ 1 - P + \frac{\sigma^2}{2} \Delta (\Delta P)' \right].
\end{eqnarray*}
The optimal $\Delta$ we seek will satisfy the two equations
\begin{eqnarray}
\frac{\partial \mathcal{L}}{\partial \Delta} - \frac{d}{d\theta} \frac{\partial \mathcal{L}}{\partial \Delta'} + \frac{d^2}{d\theta^2}\frac{\partial L}{\partial \Delta''} &=& 0 \label{eqn:EL1} \\
\frac{\partial \mathcal{L}}{\partial P} - \frac{d}{d\theta}\frac{\partial \mathcal{L}}{\partial P'} &=& 0. \label{eqn:EL2}
\end{eqnarray}
Note that we can write two more Euler-Lagrange equations, but $\frac{\partial \mathcal{L}}{\partial \nu_1} = 0$ simply restates Eq.(\ref{eqn:constraint}), and $\frac{\partial \mathcal{L}}{\partial \nu_2} = 0$ returns Eq.(\ref{eqn:steadystate}) governing $P$.

Assuming the parameter $c$ is nonzero, we obtain from Eqs.(\ref{eqn:EL1}) and (\ref{eqn:EL2}) a 4th order system of ordinary differential equations:
\begin{widetext}
\begin{eqnarray}
P'' \Delta + 2(P' \Delta' + P \Delta'' + a \Delta \nu_1 - b \Delta'' \nu_1 + c \Delta^{(4)} \nu_1) + \frac{1}{2} \Delta (P' \nu_2 - P \nu_2') \sigma^2 &=& 0 \label{eqn:4th-ode1} \\
\Delta \Delta'' - \nu_2 - \frac{1}{2} \Delta (\Delta' \nu_2 + \Delta \nu_2') \sigma^2 &=& 0.
\label{eqn:4th-ode2}
\end{eqnarray}
\end{widetext}
If $c = 0$, we will have instead the 2nd order system which
obtains by setting $c = 0$ in Eq.(\ref{eqn:4th-ode1}). 
When we examine the effects of varying the constraint 
parameters in Section \ref{sec:constraints}, we will see that 
the main result remains the same in this case as well.

\section{\label{sec:perturbation} Perturbation Approximation}
Let us first consider the 4th order case where the parameter $c$ is
nonzero. 

Assuming the noise amplitude $\sigma$ is sufficiently small, we write the following expansions
\begin{eqnarray}
P (\theta) &=& P_0 (\theta) + \sigma^2 P_1 (\theta) + ... \nonumber \\
\Delta (\theta) &=& \Delta_0 (\theta) + \sigma^2 \Delta_1 (\theta) + ...\label{eqn:perturbation} \\
\nu_1  &=& \nu_{1,0} + \sigma^2 \nu_{1,1} + ... \nonumber \\
\nu_2 (\theta)  &=& \nu_{2,0} (\theta) + \sigma^2 \nu_{2,1} (\theta)+ ... \nonumber
\end{eqnarray}
Substituting these into (\ref{eqn:4th-ode1}) and (\ref{eqn:4th-ode2}) and equating like powers of $\sigma$ gives to lowest order: $P_0(\theta) = 1$, $\nu_{2,0}(\theta) = \Delta_0(\theta) \Delta_0''(\theta)$ and the fourth order homogeneous equation
\begin{equation}
a\nu_{1,0} \Delta_0 + (1-b\nu_{1,0})\Delta_0''+c\nu_{1,0}\Delta_0^{(4)} = 0.
\label{eqn:perturbation0}
\end{equation}
For convenience let us define the differential operator
\[
\mathcal{J} = a\nu_{1,0} + (1-b\nu_{1,0})\frac{\partial^2}{\partial\theta^2} + c\nu_{1,0}\frac{\partial^4}{\partial\theta^4}.
\]
Thus Eq.(\ref{eqn:perturbation0}) becomes $\mathcal{J}(\Delta_0) = 0$, and the first order correction $\Delta_1$ obeys the inhomogeneous equation
\begin{eqnarray}
\mathcal{J}(\Delta_1)&=& (\Delta_0')^3 - b\nu_{1,1}\Delta_0'' \nonumber \\
&& + \Delta_0(a\nu_{1,1} + 3\Delta_0'\Delta_0'') + c\nu_{1,1}\Delta_0^{(4)}.
\label{eqn:perturbation1}
\end{eqnarray}
Furthermore, substituting the expansions (\ref{eqn:perturbation}) into Eq.(\ref{eqn:constraint}) gives the corresponding constraints:
\begin{eqnarray}
\int_0^1 a\Delta_0^2 + b(\Delta_0')^2+c(\Delta_0'')^2 &=& 1 \label{eqn:constraint0}\\
\int_0^1 a\Delta_0\Delta_1 + b\Delta_0'\Delta_1' + c \Delta_0''\Delta_1'' &=& 0.
\label{eqn:constraint1}
\end{eqnarray}

\begin{figure*}
 \includegraphics[width=5in,height=3.75in]{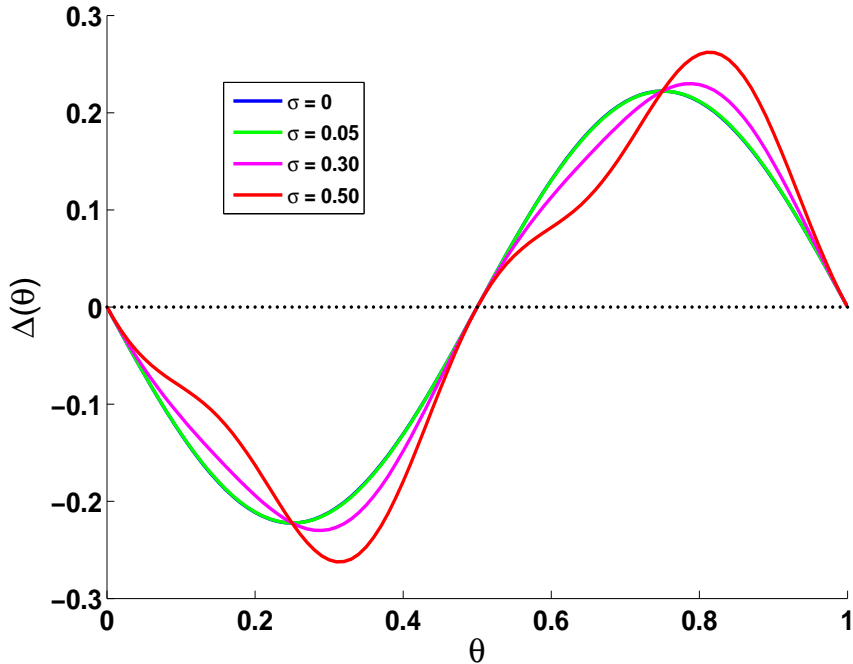}
 \caption{\label{fig:deformation} In the case where the second derivative is left unconstrained, the optimal PRC deviates from a pure cosine function as the noise amplitude $\sigma$ increases. Parameters are a=1, b=1, c=0.}
\end{figure*}

Before solving Eq.(\ref{eqn:perturbation0}), we must first determine
the unknown parameter $\nu_{1,0}$. Since we seek only periodic
solutions, we can impose a condition on the characteristic equation of
(\ref{eqn:perturbation0}): 
\begin{equation}
a\nu_{1,0} + (1-b\nu_{1,0})y^2 + c\nu_{1,0} y^4 = 0.
\label{eqn:poly}
\end{equation}
Specifically, by requiring that the roots of this polynomial satisfy $y = 2\pi i$, we determine that
\[
\nu_{1,0} = \frac{4\pi^2}{a+4b\pi^2+16c\pi^4}.
\]
Now we are ready to impose periodic boundary conditions, and we find that the solution of (\ref{eqn:perturbation0}) is just $\Delta_0(\theta) = C_0 \sin(2\pi\theta)$. The constant of integration $C_0$ is determined from the constraint (\ref{eqn:constraint0}) so that
\[
C_0 = \pm\frac{\sqrt{2}}{\sqrt{a+4b\pi^2+16c\pi^4}}.
\]
While both values of $C_0$ will give the same minimal value of the Lyapunov exponent, we choose the negative value for biological plausibility. Hence to lowest order we find the optimal phase resetting curve is Type II:
\begin{equation}
\Delta_0(\theta) = -\frac{\sqrt{2}\sin(2\pi\theta)}{\sqrt{a+4b\pi^2+16c\pi^4}}.
\label{eqn:delta0}
\end{equation}

The next order correction does not appreciably change this result. To obtain the $\sigma^2$ term, we must solve (\ref{eqn:perturbation1}) subject to (\ref{eqn:constraint1}). By the Fredholm Alternative, a solution to the inhomogeneous problem exists if and only if the right-hand side of (\ref{eqn:perturbation1}), call it $r(\theta)$, is orthogonal to the nullspace of $\mathcal{J^*}$. However, since $\mathcal{J}$ is self-adjoint we simply solve for the value of $\nu_{1,1}$ such that
\[
\int_0^1 \sin(2\pi\theta) r(\theta) d\theta = 0,
\]
namely, $\nu_{1,1} = 0$.

Imposing periodic boundary conditions on the resulting equation yields the solution
\[
\Delta_1(\theta) = C_1 \sin(2\pi\theta) + \frac{\sqrt{2}\pi\sin(2\pi\theta)\sin(4\pi\theta)}{(a-144c\pi^4)\sqrt{a+4b\pi^2+16c\pi^4}}.
\]
As before, we use the constraint (\ref{eqn:constraint1}) to obtain $C_1 = 0$. Hence to order $\sigma^2$ the optimal phase resetting curve is given by
\begin{eqnarray}
\Delta(\theta) &=& -\frac{\sqrt{2}\sin(2\pi\theta)}{\sqrt{a+4b\pi^2+16c\pi^4}} \nonumber \\
&& + \frac{\sigma^2}{2} \frac{\sqrt{2}\pi\sin(2\pi\theta)\sin(4\pi\theta)}{(a-144c\pi^4)\sqrt{a+4b\pi^2+16c\pi^4}}.
\label{eqn:optimal}
\end{eqnarray}

\section{\label{sec:constraints} Constraint Parameters}

\begin{figure*}
 \centering
 \includegraphics[width=5in,height=3.75in]{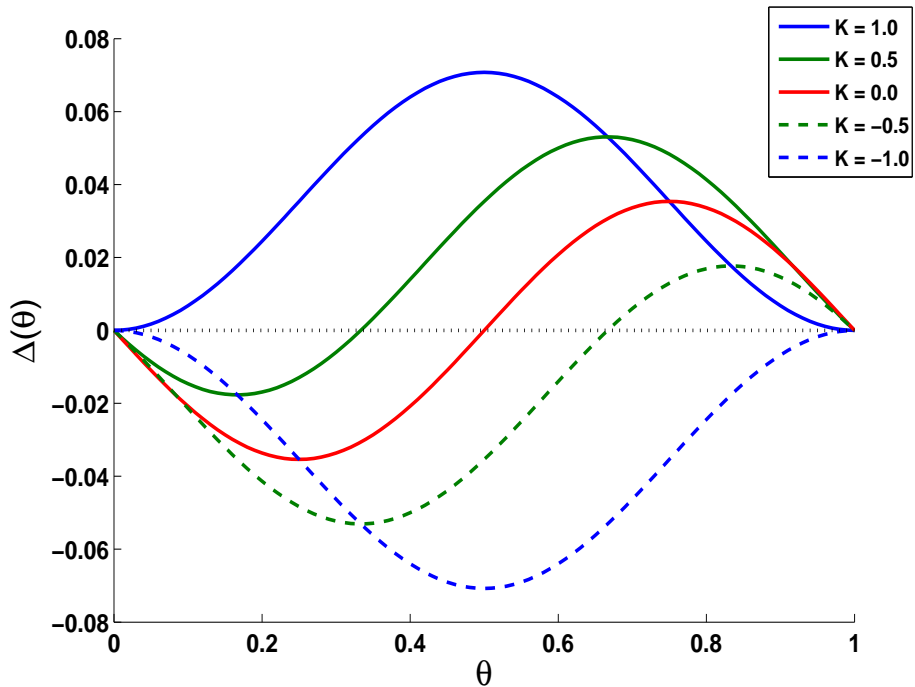}
 \caption{\label{fig:family} When the first derivative is unconstrained while the second derivative is constrained, Euler-Lagrange optimization produces a family of candidates for the minimizer of the Lyapunov exponent ranging smoothly from Type II to Type I as the parameter $K$ ranges from 0 to 1. For negative $K$ (dashed), the curves do not represent biologically plausible PRCs. Parameters are $a=0,b=1,c=1$.}
\end{figure*}

Let us next explore the influence of the constraint parameters $a$,
$b$ and $c$, which we will allow to take on the values of $0$ or
$1$. Of the seven nontrivial combinations, one has no periodic
solution at all and is thus inadmissible. Four parameter choices give
rise to the same optimum already found in Eq.(\ref{eqn:optimal}), and
two parameter combinations do not produce a unique solution but
instead yield a family of solutions ranging smoothly from Type I to
Type II. In this case, we explicitly find the minimizer of $\lambda$
among the family of solutions. 

All of the cases can be analyzed by examining Eq.(\ref{eqn:poly}), the
characteristic equation of $\mathcal{L}(\Delta)=0$. For example, the
case $a=c=0$ and $b=1$ can have no periodic solution, since the
polynomial $(1-\nu_{1,0})y^2 = 0$ has no nontrivial roots. 

The four parameter combinations that lead to Eq.(\ref{eqn:optimal})
are those in which $a = 1$. In these cases we have 
\[
\nu_{1,0} + (1-b\nu_{1,0})y^2 + c\nu_{1,0} y^4 = 0.
\]
If $c \neq 0$, the polynomial is 4th degree having four distinct roots; if $c = 0$ the polynomial is quadratic with two distinct roots. In each case we can set $y = 2\pi i$ and solve uniquely for $\nu_{1,0}$ as discussed above.

The case $c=0$ (while $a=1$) deserves further attention for
another reason. In this regime, the optimal PRC becomes sensitive to
the noise amplitude $\sigma$ as illustrated in
Fig.(\ref{fig:deformation}). To understand why the curve deforms, let
us focus on the extrema of Eq.(\ref{eqn:optimal}), which are given by
the zeros of the derivative: 
\begin{eqnarray*}
\Delta'(\theta) &=& -\frac{2\sqrt{2}\pi}{\sqrt{a+4b\pi^2+16c\pi^4}} \Bigg[\cos(2\pi \theta) \\
&& + \frac{\sigma^2\pi}{a-144c\pi^4}\Bigg(\cos(4\pi\theta)\sin(2\pi\theta) \\
&& + \frac{1}{2}\cos(2\pi\theta)\sin(4\pi\theta)\Bigg)\Bigg]
\end{eqnarray*}
In this form we clearly see that the unperturbed extrema (when $\sigma = 0$) occur at
$\theta=1/4$ and $3/4$, while deformation due to noise is on the order
of $\sigma^2\pi/(a-144c\pi^4)$. More specifically, when $c\neq0$ this
quantity is $\mathcal{O}(\sigma^2 10^{-4})$ so that the weak noise in our model
($\sigma\ll 1$) has negligible effect. However when $c=0$, this
quantity is $\mathcal{O}(\sigma^2)$, so that even relatively small magnitude
noise can have a noticible impact on the shape of the optimal PRC. 

Another interesting situation arises in the two cases where $a = 0$,
$c = 1$ and $b$ is arbitrary. Here the characteristic equation has a
double root at $y=0$: 
\[
(1-b\nu_{1,0})y^2 + \nu_{1,0} y^4 = 0.
\]
After accounting for the boundary conditions, we have a superposition of two independent solutions
\[
\Delta_0(\theta) = C_3(1-\cos(2\pi \theta)) + C_4\sin(2\pi \theta).
\]
The constraint (\ref{eqn:constraint0}) eliminates only one degree of freedom, leaving a family of solutions as candidates for the optimum:
\begin{eqnarray}
\Delta_0(\theta) &=& K \frac{1-\cos(2\pi \theta)}{\sqrt{2\pi^2 (b+4\pi^2)}} \nonumber \\
&& - \sqrt{1-K^2}\frac{\sin(2\pi \theta)}{\sqrt{2\pi^2 (b+4\pi^2)}},
\label{eqn:family}
\end{eqnarray}
where the remaining degree of freedom $K$ has been normalized to range between $-1$ and $1$. See Fig.(\ref{fig:family}).

Combining Eq.(\ref{eqn:lambda}) for the Lyapunov exponent with Eq.(\ref{eqn:ptosigma^4}) for the steady-state phase distribution, we insert Eq.(\ref{eqn:family}) to obtain the following expression:
\[
\lambda = -\frac{1}{b+4\pi^2} + \frac{\sigma^4}{4} \frac{(4K^4+10K^2+1)}{4\pi^2(b+4\pi^2)^3},
\]
where we have set $a=0,c=1$. Note that we needed to carry out the expansion of $\lambda$ to $\sigma^4$ in order to discover the dependence on $K$.

Since the derivative of $\lambda$ with respect to $K$ has only one real root at $K=0$, where a minimum occurs, the Type II curve remains the optimal PRC even in this case.

\begin{figure*}
 \centering
 \includegraphics[width=5in,height=3.75in]{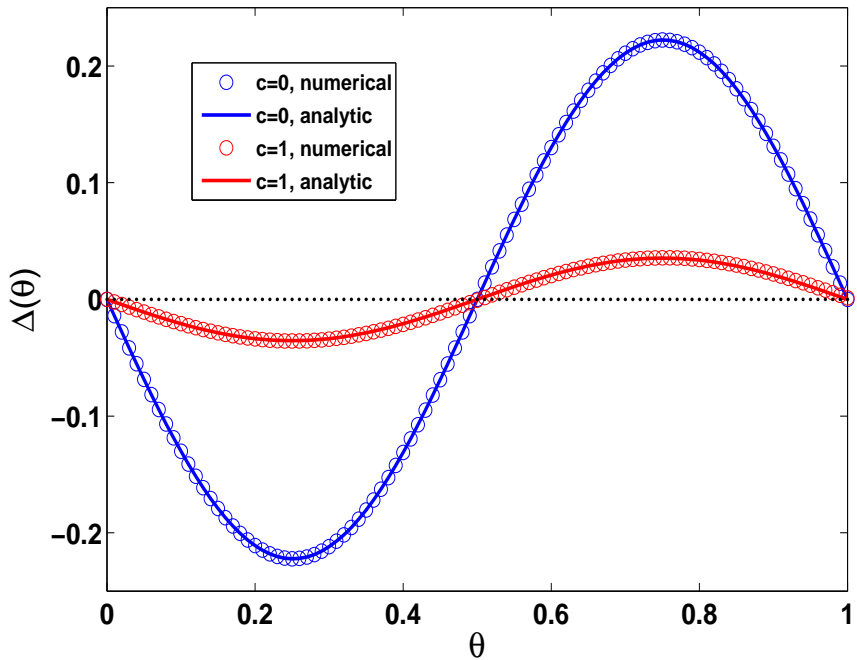}
 \caption{\label{fig:abc} The magnitude of the optimal PRC depends on the whether or not the second derivative is constrained. The numerical solution (open circles) and the analytic result (solid lines) coincide. Parameters are $a=1,b=1$ and $\sigma = 0.05$.}
\end{figure*}

\section{\label{sec:numerics} Numerical Verification}

We would like to independently verify the accuracy of the optimal PRC (\ref{eqn:optimal}) 
derived via perturbation expansion by
numerically solving the Euler-Lagrange equations (\ref{eqn:4th-ode1})
and (\ref{eqn:4th-ode2}) with periodic boundary
conditions. Unfortunately, the resulting system is singular and
therefore very difficult to solve numerically. Instead we substitute the
approximation $P(\theta) = 1 + \frac{\sigma^2}{2}\Delta(\theta)
\Delta'(\theta)$ into the Euler-Lagrange functional
(\ref{eqn:EL-functional}) to obtain a new functional 
\begin{eqnarray*}
&\int_0^1& \Delta'' \Delta \left(1 + \frac{\sigma^2}{2}\Delta(\theta) \Delta'(\theta)\right) \\
&& + \nu_1 \left[ a\Delta^2 + b(\Delta')^2 +c(\Delta'')^2 -1\right] d\theta = 0,
\label{eqn:num-functional}
\end{eqnarray*}
which gives rise via Eq.(\ref{eqn:EL1}) to the 4th order boundary value problem
\[
\Delta^{(4)} = \frac{-2\Delta''-2a\Delta\nu_1 + 2b\Delta'' \nu_1 - \Delta'^3 \sigma^2 - 3 \Delta\Delta'\Delta'' \sigma^2}{2 c \nu_1}.
\]
When $c=0$, we similarly obtain a 2nd order boundary value problem. 

Using the numerical integration package XPPAUT, we are able to achieve
excellent agreement with our analytical
approximation. Fig.(\ref{fig:abc}) illustrates numerical and analytic
solutions in the case where $c=1$ and where $c=0$. Note that imposing
a constraint on the second derivative of $\Delta$ results in an
optimal PRC of much smaller magnitude. 

\begin{figure*}
 \centering
 \includegraphics[width=5in,height=3.75in]{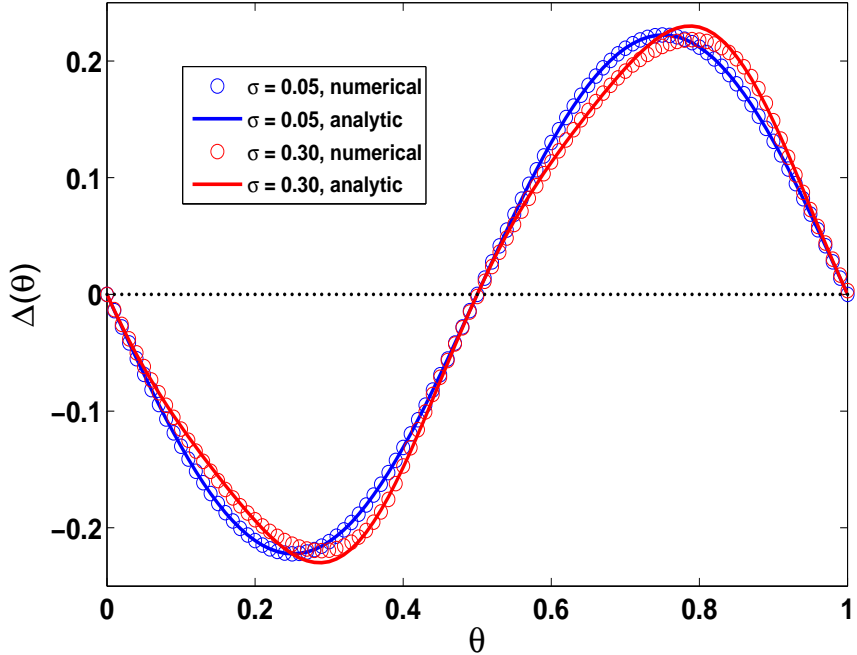}
 \caption{\label{fig:ab0} When the second derivative is unconstrained, the optimal PRC shape deforms with increasing noise. The numerical solution (open circles) and the analytic result (solid lines) are in good agreement. Parameters are $a=1,b=1,c=0$.}
\end{figure*}

In Fig.(\ref{fig:ab0}) we find good agreement between the analytic
and numerical results even for the regime in which $a=1$, $c=0$ and PRC
shape is sensitive to noise amplitude. The numerical simulation
deforms with increasing $\sigma$ just as the analytic
approximation does. 

\section*{Discussion}

In this paper we have used perturbation theory and the calculus of
variations to analyze the rate at which neurons can synchronize when
subjected to common inputs. We treat the inputs as ``noise,'' that
is, as if they are delta-correlated with no structure. Real neuronal
inputs do have correlational structure, however, so that the expression for the rate of
synchronization (the Lyapunov exponent) is more complex. 
Indeed, in previous work \cite{galanfp} we have shown that
the temporal characteristics of the noise can also have an effect on
how rapidly neurons synchronize. In that work, we asked the reverse
question: given a particular PRC, what correlation time for the noise
minimizes the Lyapunov exponent?

Suppose that we use some signal that
is not white noise but still has zero mean and is stationary. Then the
phase satisfies
\[
\frac{d\theta}{dt} = 1 + \Delta(\theta)\xi(t)
\]
where $\xi(t)$ is the input. The Lyapunov exponent is
\[
\lambda:= \lim_{T\to\infty} \frac{1}{T}\int_0^T
\Delta'(\theta(t))\xi(t)\ dt.
\] 
By using an approximation of $\theta(t)$ as in \cite{galan-optima} we
may be able to obtain a functional for $\lambda$ depending on $\xi(t)$ and
$\Delta$, and from this apply similar methods to estimate the optimal
shape of the PRC given the statistics of the inputs.

Optimization has been applied to other aspects of neural
oscillators. Moehlis, et al. \cite{moehlis} asked the following
question. Consider the scalar oscillator model:
\[
\frac{d\theta}{dt} = f(\theta) + \Delta(\theta)I(t).
\]
(Note that if $f(\theta)=1$, we have Eq.(\ref{eqn:Stratonovich}), the case considered in this paper.) Suppose the neuron fired at $t=0$ and we
desire it to fire again at time $T>0$. What is the minimum stimulus,
$I(t)$ (which, say, minimizes $\int_0^T I(t)^2 dt$) to do this? Moehlis, et al. \cite{moehlis}
write the Euler-Lagrange equations for this optimization problem and
then assume that $I(t)$ is small in order to use perturbation
methods. A related issue is the ``optimal stimulus'' \cite{rieke} for
producing a spike in a neuron, and for neural
oscillators this has been answered in \cite{ermentrout}. 


%




\appendix*
\section{An existence proof}

On the interval $[0,1]$, the phase resetting curve $\Delta$ is necessarily $0$ at the endpoints and possibly at interior points as well. As a result, we have a singular differential equation for the steady state distribution of phases $P$, derived earlier as Eq.(\ref{eqn:steadystate}) and repeated here:
\begin{equation}
 -J = -P+\frac{\sigma^2}{2}\Delta(\Delta P)'.
\label{eqn:ode4P}
\end{equation}
However we will now see that Eq.(\ref{eqn:ode4P}) does indeed have a solution despite the singularities.

Suppose $\Delta(\theta)\neq 0$ in the open interval $(a,b)\subseteq [0,1]$, while $\Delta(a)=\Delta(b)=0$. In this way, we will be able to apply our proof to the entire domain $[0,1]$ in a piecewise fashion; for example, if $\Delta(x) = \sin(2\pi x)$, then $a=0$ and $b=1/2$, or $a=1/2$ and $b=1$. In the following we will assume, without loss of generality, that $\Delta(\theta) > 0$ in $(a,b)$.

Let us begin by rewriting the differential equation as an integral equation. Define $Q(x) := \Delta(x) P(x).$ Then Eq.(\ref{eqn:ode4P}) becomes
\begin{equation}
Q' - \frac{2 Q}{\sigma^2 \Delta^2} = \frac{-2 J}{\sigma^2 \Delta}.
\label{eqn:ode4Q}
\end{equation}
We now introduce an integrating factor; let
\[
z(x) := -\frac{2}{\sigma^2}\int_c^x \frac{ds}{\Delta^2(s)},
\]
where $c\in (a,b)$ is fixed. Observe that, as $x$ approaches $a$ from above we eventually have $x<c$, and hence $z(x)$ approaches $+\infty$. Likewise, as $x$ approaches $b$ from below, $z(x)$ approaches $-\infty$.

Eq.(\ref{eqn:ode4Q}) now becomes
\[
(e^{z(x)} Q)' = -\frac{2 J}{\sigma^2 \Delta} e^{z(x)}.
\]
Integrating both sides gives
\begin{equation}
Q(x) = \frac{2 J}{\sigma^2} e^{-z(x)}\left( K-\int_c^x \frac{e^{z(t)}}{\Delta(t)} dt \right),
\label{eqn:integral}
\end{equation}
where $K$ is a constant of integration that will be determined below.

We see from Eq.(\ref{eqn:ode4P}) that $P(a)=P(b)=J$. Therefore a solution exists iff $\lim_{x\rightarrow a^+} Q(x)/\Delta(x) = \lim_{x\rightarrow b^-} Q(x)/\Delta(x) = J$.
Let us first consider the right endpoint and assume for now that the limit
\begin{equation}
\lim_{x\rightarrow b^-} \int_c^x \frac{e^{z(t)}}{\Delta(t)} dt = L
\label{eqn:claim}
\end{equation}
exists. Let us compute 
\begin{eqnarray*}
\lim_{x\rightarrow b^-}\frac{Q(x)}{\Delta(x)} &=& \frac{2 J}{\sigma^2} \lim_{x\rightarrow b^-} \frac{e^{-z(x)}}{\Delta(x)} \left( K-\int_c^x \frac{e^{z(t)}}{\Delta(t)} dt \right) \\
 &=& \frac{2 J}{\sigma^2} \lim_{x\rightarrow b^-} \frac{K-\int_c^x \frac{e^{z(t)}}{\Delta(t)} dt}{\Delta(x)e^{z(x)}},
\end{eqnarray*}
and note that when we set $K = L$, both numerator and denominator tend to 0 as $x\rightarrow b^-$. Thus we can use L'H\^{o}pital's rule to obtain
\begin{eqnarray}
\lim_{x\rightarrow b^-} \frac{Q(x)}{\Delta(x)} &=& \frac{2 J}{\sigma^2} \lim_{x\rightarrow b^-} \frac{-e^{z(x)}/\Delta(x)}{\Delta(x)z'(x) e^{z(x)}+\Delta'(x)e^{z(x)}}  \nonumber \\
 &=& \frac{2 J}{\sigma^2} \lim_{x\rightarrow b^-} \frac{-1/\Delta(x)}{-\frac{2}{\sigma^2}\frac{1}{\Delta(x)^2}\Delta(x) +\Delta'(x) } \nonumber \\
 &=& \frac{2 J}{\sigma^2} \lim_{x\rightarrow b^-} \frac{-1}{-\frac{2}{\sigma^2} +\Delta(x)\Delta'(x) } \nonumber \\ 
 &=& \frac{2 J}{\sigma^2} \left(\frac{\sigma^2}{2}\right) \nonumber \\ 
 &=& J.
\label{eqn:lhopitals}
\end{eqnarray}

Now let us return to the assumption we made and observe that the integral in Eq.(\ref{eqn:claim}) is not improper after all. Rewriting the integrand of (\ref{eqn:claim}) such that both numerator and denominator go to infinity, we can use L'H\^{o}pital's rule again to see that the integrand goes to zero:
\begin{eqnarray*}
\lim_{t\rightarrow b^-} \frac{e^{z(t)}}{\Delta(t)} 
 &=& \lim_{t\rightarrow b^-} \frac{1/\Delta(t)}{e^{-z(t)}} \\
 &=& \lim_{t\rightarrow b^-} \frac{-\Delta'(t)/\Delta(t)^2}{e^{-z(t)}/\Delta(t)^2} \\
 &=& \lim_{t\rightarrow b^-} -\Delta'(t) e^{z(t)}\\
 &=& 0.
\end{eqnarray*}
The last equality follows since $\Delta'$ is bounded and $\lim_{x\rightarrow b^-} e^{z(t)} = 0$. Hence our assumption was justified.

Now let us rewrite Eq.(\ref{eqn:integral}), incorporating our knowledge from Eq.(\ref{eqn:claim}), namely that $K=L$:
\begin{eqnarray*}
Q(x) &=& \frac{2 J}{\sigma^2} e^{-z(x)}\left( \int_c^b \frac{e^{z(t)}}{\Delta(t)} dt - \int_c^x \frac{e^{z(t)}}{\Delta(t)} dt \right) \nonumber \\
 &=& \frac{2 J}{\sigma^2} e^{-z(x)} \int_x^b \frac{e^{z(t)}}{\Delta(t)} dt.
\end{eqnarray*}
It remains to show that $\lim_{x\rightarrow a^+} Q(x)/\Delta(x) = J$. We will prepare to use L'H\^{o}pital's rule once again by writing
\begin{eqnarray}
\lim_{x\rightarrow a^+} \frac{Q(x)}{\Delta(x)} &=& \frac{2 J}{\sigma^2} \lim_{x\rightarrow a^+} \frac{e^{-z(x)}}{\Delta(x)} \int_x^b \frac{e^{z(t)}}{\Delta(t)} dt \nonumber \\
 &=& \frac{2 J}{\sigma^2} \lim_{x\rightarrow a^+} \frac{\int_x^b \frac{e^{z(t)}}{\Delta(t)}dt}{\Delta(x) e^{z(x)}}.
\label{eqn:xapproachesa}
\end{eqnarray}
Since $e^{z(t)}$ tends to infinity as $x$ approaches $a$ from above, by L'H\^{o}pital's rule the denominator of (\ref{eqn:xapproachesa}) also tends to infinity:
\begin{eqnarray*}
\lim_{x\rightarrow a^+} \frac{e^{z(x)}}{1/\Delta(x)} &=& -\frac{2}{\sigma^2}\lim_{x\rightarrow a^+} \frac{e^{z(x)}/\Delta(x)^2}{\Delta'(x)/\Delta(x)^2} \\
 &=& -\frac{2}{\sigma^2}\lim_{x\rightarrow a^+} \frac{e^{z(x)}}{\Delta'(x)} \\
 &=& \infty.
\end{eqnarray*}
The numerator of Eq.(\ref{eqn:xapproachesa}) tends to infinity as well since
\[
\int_x^b \frac{e^{z(t)}}{\Delta(t)} dt > \int_x^b \frac{e^{z(t)}}{M} dt,
\]
when $M=\max\{\Delta(x):x\in [0,1]\}$, and the latter integral is clearly unbounded as $x$ approaches $a$. Therefore we can apply to (\ref{eqn:xapproachesa}) a similar calculation to that in (\ref{eqn:lhopitals}) and conclude that $\lim_{x\rightarrow a^+} Q(x)/\Delta(x) = J$ as desired.


\bibliography{optimalprc_aps}

\begin{thebibliography}{23}
\expandafter\ifx\csname natexlab\endcsname\relax\def\natexlab#1{#1}\fi
\expandafter\ifx\csname bibnamefont\endcsname\relax
  \def\bibnamefont#1{#1}\fi
\expandafter\ifx\csname bibfnamefont\endcsname\relax
  \def\bibfnamefont#1{#1}\fi
\expandafter\ifx\csname citenamefont\endcsname\relax
  \def\citenamefont#1{#1}\fi
\expandafter\ifx\csname url\endcsname\relax
  \def\url#1{\texttt{#1}}\fi
\expandafter\ifx\csname urlprefix\endcsname\relax\def\urlprefix{URL }\fi
\providecommand{\bibinfo}[2]{#2}
\providecommand{\eprint}[2][]{\url{#2}}

\bibitem[{\citenamefont{Tiesinga}(2004)}]{tiesinga:2004}
\bibinfo{author}{\bibfnamefont{P.~H.~E.} \bibnamefont{Tiesinga}},
  \bibinfo{journal}{Phys. Rev. E} \textbf{\bibinfo{volume}{69}},
  \bibinfo{pages}{031912} (\bibinfo{year}{2004}).

\bibitem[{\citenamefont{Teramae and Tanaka}(2004)}]{teramae:2004}
\bibinfo{author}{\bibfnamefont{J.~N.} \bibnamefont{Teramae}} \bibnamefont{and}
  \bibinfo{author}{\bibfnamefont{D.}~\bibnamefont{Tanaka}},
  \bibinfo{journal}{Phys. Rev. Lett.} \textbf{\bibinfo{volume}{93}},
  \bibinfo{pages}{204103} (\bibinfo{year}{2004}).

\bibitem[{\citenamefont{Goldobin and Pikovsky}(2005)}]{goldobin:2005}
\bibinfo{author}{\bibfnamefont{D.~S.} \bibnamefont{Goldobin}} \bibnamefont{and}
  \bibinfo{author}{\bibfnamefont{A.}~\bibnamefont{Pikovsky}},
  \bibinfo{journal}{Phys. Rev. E} \textbf{\bibinfo{volume}{71}},
  \bibinfo{pages}{045201(R)} (\bibinfo{year}{2005}).

\bibitem[{\citenamefont{Nakao et~al.}(2005)\citenamefont{Nakao, Arai, Nagai,
  Tsubo, and Kuramoto}}]{nakao:2005}
\bibinfo{author}{\bibfnamefont{H.}~\bibnamefont{Nakao}},
  \bibinfo{author}{\bibfnamefont{K.~S.} \bibnamefont{Arai}},
  \bibinfo{author}{\bibfnamefont{K.}~\bibnamefont{Nagai}},
  \bibinfo{author}{\bibfnamefont{Y.}~\bibnamefont{Tsubo}}, \bibnamefont{and}
  \bibinfo{author}{\bibfnamefont{Y.}~\bibnamefont{Kuramoto}},
  \bibinfo{journal}{Phys. Rev. E} \textbf{\bibinfo{volume}{72}},
  \bibinfo{pages}{026220} (\bibinfo{year}{2005}).

\bibitem[{\citenamefont{Stroeve and Gielen}(2001)}]{stroeve:2001}
\bibinfo{author}{\bibfnamefont{S.}~\bibnamefont{Stroeve}} \bibnamefont{and}
  \bibinfo{author}{\bibfnamefont{S.}~\bibnamefont{Gielen}},
  \bibinfo{journal}{Neural Comput.} \textbf{\bibinfo{volume}{13}},
  \bibinfo{pages}{2005} (\bibinfo{year}{2001}).

\bibitem[{\citenamefont{Engel et~al.}(1991{\natexlab{a}})\citenamefont{Engel,
  Kreiter, Konig, and Singer}}]{engel:1991}
\bibinfo{author}{\bibfnamefont{A.~K.} \bibnamefont{Engel}},
  \bibinfo{author}{\bibfnamefont{A.~K.} \bibnamefont{Kreiter}},
  \bibinfo{author}{\bibfnamefont{P.}~\bibnamefont{Konig}}, \bibnamefont{and}
  \bibinfo{author}{\bibfnamefont{W.}~\bibnamefont{Singer}},
  \bibinfo{journal}{Proc. Natl. Acad. Sci.} \textbf{\bibinfo{volume}{88}},
  \bibinfo{pages}{6048} (\bibinfo{year}{1991}{\natexlab{a}}).

\bibitem[{\citenamefont{Engel et~al.}(1991{\natexlab{b}})\citenamefont{Engel,
  Konig, Kreiter, and Singer}}]{engel:1991b}
\bibinfo{author}{\bibfnamefont{A.~K.} \bibnamefont{Engel}},
  \bibinfo{author}{\bibfnamefont{P.}~\bibnamefont{Konig}},
  \bibinfo{author}{\bibfnamefont{A.~K.} \bibnamefont{Kreiter}},
  \bibnamefont{and} \bibinfo{author}{\bibfnamefont{W.}~\bibnamefont{Singer}},
  \bibinfo{journal}{Science} \textbf{\bibinfo{volume}{252}},
  \bibinfo{pages}{1177} (\bibinfo{year}{1991}{\natexlab{b}}).

\bibitem[{\citenamefont{Gal\'{a}n et~al.}(2006)\citenamefont{Gal\'{a}n,
  Fourcaud-Trocme, Ermentrout, and Urban}}]{galan:2006}
\bibinfo{author}{\bibfnamefont{R.~F.} \bibnamefont{Gal\'{a}n}},
  \bibinfo{author}{\bibfnamefont{N.}~\bibnamefont{Fourcaud-Trocme}},
  \bibinfo{author}{\bibfnamefont{G.~B.} \bibnamefont{Ermentrout}},
  \bibnamefont{and} \bibinfo{author}{\bibfnamefont{N.~N.} \bibnamefont{Urban}},
  \bibinfo{journal}{J. Neurosci.} \textbf{\bibinfo{volume}{26}},
  \bibinfo{pages}{3646} (\bibinfo{year}{2006}).

\bibitem[{\citenamefont{Hansel et~al.}(1995)\citenamefont{Hansel, Mato, and
  Meunier}}]{hansel:1995}
\bibinfo{author}{\bibfnamefont{D.}~\bibnamefont{Hansel}},
  \bibinfo{author}{\bibfnamefont{G.}~\bibnamefont{Mato}}, \bibnamefont{and}
  \bibinfo{author}{\bibfnamefont{C.}~\bibnamefont{Meunier}},
  \bibinfo{journal}{Neural Comput.} \textbf{\bibinfo{volume}{7}},
  \bibinfo{pages}{307} (\bibinfo{year}{1995}).

\bibitem[{\citenamefont{Ermentrout et~al.}(2001)\citenamefont{Ermentrout,
  Pascal, and Gutkin}}]{ermentrout:2001}
\bibinfo{author}{\bibfnamefont{G.~B.} \bibnamefont{Ermentrout}},
  \bibinfo{author}{\bibfnamefont{M.}~\bibnamefont{Pascal}}, \bibnamefont{and}
  \bibinfo{author}{\bibfnamefont{B.~S.} \bibnamefont{Gutkin}},
  \bibinfo{journal}{Neural Comput.} \textbf{\bibinfo{volume}{13}},
  \bibinfo{pages}{1285} (\bibinfo{year}{2001}).

\bibitem[{\citenamefont{Gutkin et~al.}(2005)\citenamefont{Gutkin, Ermentrout,
  and Reyes}}]{gutkin:2005}
\bibinfo{author}{\bibfnamefont{B.~S.} \bibnamefont{Gutkin}},
  \bibinfo{author}{\bibfnamefont{G.~B.} \bibnamefont{Ermentrout}},
  \bibnamefont{and} \bibinfo{author}{\bibfnamefont{A.~D.} \bibnamefont{Reyes}},
  \bibinfo{journal}{J. Neurophysiol.} \textbf{\bibinfo{volume}{94}},
  \bibinfo{pages}{1623} (\bibinfo{year}{2005}).

\bibitem[{\citenamefont{Netoff et~al.}(2005)\citenamefont{Netoff, Acker,
  Bettencourt, and White}}]{netoff:2005}
\bibinfo{author}{\bibfnamefont{T.~I.} \bibnamefont{Netoff}},
  \bibinfo{author}{\bibfnamefont{C.~D.} \bibnamefont{Acker}},
  \bibinfo{author}{\bibfnamefont{J.~C.} \bibnamefont{Bettencourt}},
  \bibnamefont{and} \bibinfo{author}{\bibfnamefont{J.~A.} \bibnamefont{White}},
  \bibinfo{journal}{J. Comput. Neurosci.} \textbf{\bibinfo{volume}{18}},
  \bibinfo{pages}{287} (\bibinfo{year}{2005}).

\bibitem[{\citenamefont{Izhikevich}(2006)}]{izhikevich}
\bibinfo{author}{\bibfnamefont{E.~M.} \bibnamefont{Izhikevich}},
  \emph{\bibinfo{title}{Dynamical Systems in Neuroscience: The Geometry of
  Excitability and Bursting}} (\bibinfo{publisher}{MIT Press},
  \bibinfo{year}{2006}).

\bibitem[{\citenamefont{Brown et~al.}(2004)\citenamefont{Brown, Moehlis, and
  Holmes}}]{brown}
\bibinfo{author}{\bibfnamefont{E.}~\bibnamefont{Brown}},
  \bibinfo{author}{\bibfnamefont{J.}~\bibnamefont{Moehlis}}, \bibnamefont{and}
  \bibinfo{author}{\bibfnamefont{P.}~\bibnamefont{Holmes}},
  \bibinfo{journal}{Neural Comp.} \textbf{\bibinfo{volume}{16}},
  \bibinfo{pages}{673} (\bibinfo{year}{2004}).

\bibitem[{\citenamefont{Gal\'{a}n et~al.}(2007)\citenamefont{Gal\'{a}n,
  Ermentrout, and Urban}}]{galanfp}
\bibinfo{author}{\bibfnamefont{R.~F.} \bibnamefont{Gal\'{a}n}},
  \bibinfo{author}{\bibfnamefont{G.~B.} \bibnamefont{Ermentrout}},
  \bibnamefont{and} \bibinfo{author}{\bibfnamefont{N.~N.} \bibnamefont{Urban}},
  \bibinfo{journal}{Phys. Rev. E} \textbf{\bibinfo{volume}{76}},
  \bibinfo{pages}{056110} (\bibinfo{year}{2007}).

\bibitem[{\citenamefont{Marella and Ermentrout}(2008)}]{sashi}
\bibinfo{author}{\bibfnamefont{S.}~\bibnamefont{Marella}} \bibnamefont{and}
  \bibinfo{author}{\bibfnamefont{G.~B.} \bibnamefont{Ermentrout}},
  \bibinfo{journal}{Phys. Rev. E} \textbf{\bibinfo{volume}{77}},
  \bibinfo{pages}{041918} (\bibinfo{year}{2008}).

\bibitem[{\citenamefont{Kuramoto}(1984)}]{kuramoto:1984}
\bibinfo{author}{\bibfnamefont{Y.}~\bibnamefont{Kuramoto}},
  \emph{\bibinfo{title}{Chemical Oscillation, Waves and Turbulence}}
  (\bibinfo{publisher}{Springer-Verlag}, \bibinfo{year}{1984}).

\bibitem[{\citenamefont{Horsthemke and Lefever}(1984)}]{horsthemke:1984}
\bibinfo{author}{\bibfnamefont{W.}~\bibnamefont{Horsthemke}} \bibnamefont{and}
  \bibinfo{author}{\bibfnamefont{R.}~\bibnamefont{Lefever}},
  \emph{\bibinfo{title}{Noise-Induced Transitions}}
  (\bibinfo{publisher}{Springer-Verlag}, \bibinfo{year}{1984}).

\bibitem[{\citenamefont{Yoshimura and Arai}(2008)}]{yoshi}
\bibinfo{author}{\bibfnamefont{K.}~\bibnamefont{Yoshimura}} \bibnamefont{and}
  \bibinfo{author}{\bibfnamefont{K.}~\bibnamefont{Arai}},
  \bibinfo{journal}{Phys. Rev. Lett.} \textbf{\bibinfo{volume}{101}},
  \bibinfo{pages}{154101} (\bibinfo{year}{2008}).

\bibitem[{\citenamefont{Gal\'{a}n et~al.}(2008)\citenamefont{Gal\'{a}n,
  Ermentrout, and Urban}}]{galan-optima}
\bibinfo{author}{\bibfnamefont{R.~F.} \bibnamefont{Gal\'{a}n}},
  \bibinfo{author}{\bibfnamefont{G.~B.} \bibnamefont{Ermentrout}},
  \bibnamefont{and} \bibinfo{author}{\bibfnamefont{N.~N.} \bibnamefont{Urban}},
  \bibinfo{journal}{J. Neurophysiol.} \textbf{\bibinfo{volume}{99}},
  \bibinfo{pages}{277} (\bibinfo{year}{2008}).

\bibitem[{\citenamefont{Moehlis et~al.}(2006)\citenamefont{Moehlis, Shea-Brown,
  and Rabitz}}]{moehlis}
\bibinfo{author}{\bibfnamefont{J.}~\bibnamefont{Moehlis}},
  \bibinfo{author}{\bibfnamefont{E.}~\bibnamefont{Shea-Brown}},
  \bibnamefont{and} \bibinfo{author}{\bibfnamefont{H.}~\bibnamefont{Rabitz}},
  \bibinfo{journal}{ASME J. of Computational and Nonlinear Dynamics}
  \textbf{\bibinfo{volume}{1}}, \bibinfo{pages}{358} (\bibinfo{year}{2006}).

\bibitem[{\citenamefont{Rieke et~al.}(1999)\citenamefont{Rieke, Warland, van
  Steveninck, and Bialek}}]{rieke}
\bibinfo{author}{\bibfnamefont{F.}~\bibnamefont{Rieke}},
  \bibinfo{author}{\bibfnamefont{D.}~\bibnamefont{Warland}},
  \bibinfo{author}{\bibfnamefont{R.}~\bibnamefont{van Steveninck}},
  \bibnamefont{and} \bibinfo{author}{\bibfnamefont{W.}~\bibnamefont{Bialek}},
  \emph{\bibinfo{title}{Spikes: Exploring the Neural Code}}
  (\bibinfo{publisher}{MIT Press}, \bibinfo{year}{1999}).

\bibitem[{\citenamefont{Ermentrout et~al.}(2007)\citenamefont{Ermentrout,
  Gal\'{a}n, and Urban}}]{ermentrout}
\bibinfo{author}{\bibfnamefont{G.~B.} \bibnamefont{Ermentrout}},
  \bibinfo{author}{\bibfnamefont{R.~F.} \bibnamefont{Gal\'{a}n}},
  \bibnamefont{and} \bibinfo{author}{\bibfnamefont{N.~N.} \bibnamefont{Urban}},
  \bibinfo{journal}{Phys. Rev. Lett.} \textbf{\bibinfo{volume}{99}},
  \bibinfo{pages}{248103} (\bibinfo{year}{2007}).

\end{thebibliography}

\end{document}